\documentclass[reqno]{amsart}
\usepackage{hyperref}

\AtBeginDocument{{\noindent\small
\emph{Electronic Journal of Differential Equations},
Vol. 2015 (2015), No. 27, pp. 1--8.\newline
ISSN: 1072-6691. URL: http://ejde.math.txstate.edu or http://ejde.math.unt.edu
\newline ftp ejde.math.txstate.edu}
\thanks{\copyright 2015 Texas State University - San Marcos.}
\vspace{9mm}}

\begin{document}
\title[\hfilneg EJDE-2015/27\hfil Inverse Sturm-Liouville problems]
{Inverse Sturm-Liouville problems with fixed boundary conditions}

\author[Yu. A. Ashrafyan,   T. N. Harutyunyan \hfil EJDE-2015/27\hfilneg]
{Yuri A. Ashrafyan,  Tigran N. Harutyunyan}

\address{Yuri A. Ashrafyan \newline
Yerevan State University, Armenia}
\email{yuriashrafyan@ysu.am}

\address{Tigran N. Harutyunyan \newline
Yerevan State University, Armenia}
\email{hartigr@yahoo.co.uk}

\thanks{Submitted December 17, 2014. Published January 28, 2015.}

\begin{abstract}
 Necessary and sufficient conditions for two sequences
 $\{\mu_n\}_{n=0}^\infty$ and $\{ a_n\}_{n=0}^\infty$ to be the spectral
 data for a certain Sturm-Liouville problem are well known.
 We add two more conditions so that the same two sequences become necessary
 and sufficient for being the  spectral data for a Sturm-Liouville problem
 with fixed boundary conditions.
\end{abstract}

\maketitle
\numberwithin{equation}{section}
\newtheorem{theorem}{Theorem}[section]
\allowdisplaybreaks

\section{Introduction and statements of the results}

 Let us denote by $L(q, \alpha, \beta)$ the Sturm-Liouville boundary-value problem
\begin{gather}
\ell y\equiv -y''+q(x)y=\mu y,\quad x\in (0, \pi),\; \mu\in \mathbb{C},\label{eq1}\\
y(0)\cos\alpha+y'(0)\sin\alpha=0,\quad \alpha\in (0, \pi],\label{eq2}\\
y(\pi)\cos\beta +y'(\pi)\sin\beta=0,\quad \beta\in[0, \pi),\label{eq3}
\end{gather}
where $q$ is a real-valued functions which are integrable on $[0, \pi]$
(we write $q\in L^1_{\mathbb{R}}[0, \pi]$).
By $L(q, \alpha, \beta)$ we also denote the self-adjoint operator, generated by
problem \eqref{eq1}-\eqref{eq3} (see \cite{n1}).
It is known, that under these conditions the spectra of the operator
 $L(q, \alpha, \beta)$ is discrete and consists of real, simple eigenvalues
\cite{n1}, which we denote by $\mu_n=\mu_n(q,\alpha,\beta)=\lambda_n^2(q, \alpha, \beta)$,
$n=0, 1, 2, \dots$, emphasizing the dependence of $\mu_n$ on $q$, $\alpha$
and $\beta$. We assume that eigenvalues are enumerated in the increasing order, i.e.,
$$
\mu_0(q, \alpha, \beta) < \mu_1(q, \alpha, \beta) < \dots < \mu_n(q, \alpha, \beta) < \dots.
$$

Let $\varphi(x,\mu,\alpha,q)$ and $\psi(x,\mu,\beta,q)$ be the solutions of the
equation \eqref{eq1}, which satisfy the initial conditions
\begin{gather*}
\varphi(0,\mu,\alpha,q)=\sin\alpha,\quad \varphi'(0,\mu,\alpha,q)=-\cos\alpha,\\
\psi(\pi,\mu,\beta,q)=\sin\beta, \quad  \psi'(\pi,\mu,\beta,q)=-\cos\beta,
\end{gather*}
respectively.
The eigenvalues $\mu_n=\mu_n(q, \alpha, \beta)$, $n=0,1, 2, \dots$, of
$L(q, \alpha, \beta)$ are the solutions of the equation
$$
\Phi(\mu)=\Phi(\mu,\alpha,\beta):=\varphi(\pi, \mu, \alpha)\cos \beta+\varphi'(\pi, \mu,\alpha)\sin\beta=0,
$$
or of the equation
$$
\Psi(\mu)=\Psi(\mu,\alpha,\beta):=\psi(0, \mu, \beta)\cos\alpha+\psi'(0, \mu, \beta)\sin \alpha=0.
$$
According to the well-known Liouville formula, the wronskian
$W(x)=W(x,\varphi,\psi)= \varphi \psi'-\varphi'\psi$ of the solutions $\varphi$ and
 $\psi$ is constant. It follows that $W(0)=W(\pi)$ and, consequently
$\Psi(\mu, \alpha, \beta)=-\Phi(\mu, \alpha, \beta)$.
It is easy to see that the functions $\varphi_n(x):=\varphi(x, \mu_n, \alpha, q)$ and
$\psi_n(x):=\psi(x,\mu_n, \beta, q)$, $n=0, 1, 2, \dots$, are the eigenfunctions,
corresponding to the eigenvalue $\mu_n$.
Since all eigenvalues are simple, there exist constants $c_n=c_n(q, \alpha, \beta)$,
$n=0, 1, 2, \dots$, such that
\begin{equation}\label{eq4}
\varphi_n(x)=c_n \psi_n(x).
\end{equation}
The squares of the $L^2$-norm of these eigenfunctions:
\begin{gather*}
a_n=a_n(q,\alpha,\beta):=\int_0^{\pi} |\varphi_n(x)|^2 dx,\quad    n=0, 1, 2, \dots,\\
b_n=b_n(q,\alpha,\beta):=\int_0^{\pi} |\psi_n(x)|^2 dx,\quad   n=0, 1, 2, \dots
\end{gather*}
are called norming constants.

In this article we consider the case $\alpha,\beta \in (0,\pi)$;
i.e. we assume that $\sin \alpha \neq 0$ and $\sin \beta \neq 0$.
In this case we consider the solution
$\tilde{\varphi}(x,\mu,\alpha,q):=\frac{\varphi(x,\mu,\alpha,q)}{\sin\alpha}$
of \eqref{eq1} which has the initial values
$$
\tilde{\varphi}(0,\mu,\alpha,q)=1, \quad \tilde{\varphi}'(0,\mu,\alpha,q)=-\cot\alpha;
$$
also we consider the solution
$\tilde{\psi}(x,\mu,\beta,q):=\frac{\psi(x,\mu,\beta,q)}{\sin\beta}$.
Of course, the functions $\tilde{\varphi}_n(x):=\tilde{\varphi}(x,\mu_n,\alpha,q)$ and
$\tilde{\psi}_n(x):=\tilde{\psi}(x,\mu_n,\alpha,q)$, $n=0,1,2,\dots$,
are the eigenfunctions, corresponding to the eigenvalue $\mu_n$.
It follows from \eqref{eq4} that for norming constants
$\tilde{a}_n:=\|\tilde{\varphi}_n\|^2=\frac{a_n}{\sin^2\alpha}$ and
$\tilde{b}_n:=\|\tilde{\psi}_n\|^2=\frac{b_n}{\sin^2\beta}$ satisfy
\begin{equation}\label{eq5}
\tilde{b}_n=\frac{b_n}{\sin^2\beta}=\frac{a_n}{c^2_n \sin^2\beta}
=\frac{\tilde{a}_n \sin^2 \alpha}{c^2_n \sin^2\beta}.
\end{equation}

The inverse problem by ``spectral function''
(see \cite{f1,g1,g2,i1,l1,m1,p1,z1})
is  the reconstruction of the problem $(q,\alpha,\beta)$ from the spectra
$\{\mu_n\} _{n=0}^{\infty}$ and the norming constants
$\{\tilde{a}_n \} _{n=0}^{\infty}$ (or $\{\tilde{b}_n \} _{n=0}^{\infty}$).
The two sequences
$\{\mu_n \} _{n=0}^\infty$ and $\{\tilde{a}_n \} _{n=0}^\infty$
together will be  called the spectral data.

In this article we state the question
\begin{quote}
What kind of sequences $\{\mu_n \} _{n=0}^\infty$ and
$\{\tilde{a}_n \} _{n=0}^\infty$ can be the spectral data for
 problem $(q,\alpha,\beta)$ with $q\in L^2_{\mathbb{R}}[0, \pi]$ 
and in advance fixed $\alpha$ and $\beta$ in $(0,\pi)$?
\end{quote}
Our answer is in the following assertion.

\begin{theorem}\label{thm1}
 For a real increasing sequence $\{\mu_n \} _{n=0}^\infty$ and a positive
sequence $\{\tilde{a}_n \} _{n=0}^\infty$
to be spectral data for boundary-value problem $(q,\alpha,\beta)$ with a
$q\in L^2_{\mathbb{R}}[0, \pi]$ and fixed $\alpha,\beta \in (0,\pi)$ it is necessary
and sufficient that the following relations hold:
\begin{gather}
 \lambda_n=\sqrt{\mu_n}=n+\frac{\omega}{\pi n}+\frac{\omega_n}{n}, \quad
 \omega=\textrm{const}, \quad \{\omega_n \} _{n=0}^\infty \in l^2, \label{eq6}
\\
 \tilde{a}_n=\frac{\pi}{2}+\frac{\kappa_n}{n}, \quad
 \{{\kappa_n} \}_{n=0}^\infty \in l^2, \label{eq7}
\\
\frac{1}{\tilde{a}_0}-\frac{1}{\pi}+\sum_{n=1}^\infty\Big(\frac{1}{\tilde{a}_n}
 -\frac{2}{\pi}\Big)=\cot\alpha, \label{eq8}
\\
 \frac{\tilde{a}_0}{\pi^2 \cdot  \left( \prod_{k=1}^\infty
 \frac{\mu_k - \mu_0}{k^2}\right)^2}-\frac{1}{\pi}
+\sum_{n=1}^\infty\Big(\frac{\tilde{a}_n n^4}{\pi^2  [\mu_0-\mu_n]^2
 \big( \prod_{k=1, k \neq n}^\infty \frac{\mu_k - \mu_n}{k^2}\big)^2}-
\frac{2}{\pi}\Big)=-\cot\beta.\label{eq9}
\end{gather}
\end{theorem}

To prove Theorem \ref{thm1} we  prove the following assertion,
 which has independent interest.

\begin{theorem}\label{thm2}
Let $q\in L^2_{\mathbb{R}}[0, \pi]$ and $\alpha,\beta \in (0,\pi)$.
Then for norming constants $\tilde{a}_n=\tilde{a}_n(q, \alpha,\beta)$ and
$\tilde{b}_n=\tilde{b}_n(q,\alpha,\beta)$ satisfy
\begin{gather}
\frac{1}{\tilde{a}_0}-\frac{1}{\pi}+\sum_{n=1}^\infty
\Big(\frac{1}{\tilde{a}_n}-\frac{2}{\pi}\Big)=\cot\alpha, \label{eq10}\\
\frac{1}{\tilde{b}_0}-\frac{1}{\pi}+\sum_{n=1}^\infty
\Big(\frac{1}{\tilde{b}_n}-\frac{2}{\pi}\Big)=-\cot\beta. \label{eq11}
\end{gather}
\end{theorem}

Let us note that asymptotic behavior of $\{\mu_n\}_{n=0}^\infty$ and
 $\{\tilde{a}_n\}_{n=0}^\infty$ are standard conditions for the solution of
the inverse problem.
The conditions \eqref{eq8} and \eqref{eq9} which we add to the conditions
\eqref{eq6} and \eqref{eq7} guarantee that $\alpha$ and $\beta$, which we construct
 during the solution of the inverse problem, are the same that we
fixed in advance. At the same time Theorem \ref{thm2} says that the conditions
\eqref{eq8} and \eqref{eq9} are necessary.

\section{Proof of Theorem \ref{thm2}}

The solution $\tilde{\varphi}$ has the well known representation
(see \cite{f1,g1,g2,l1,m1})
\begin{equation}\label{eq12}
\tilde{\varphi}(x,\lambda,\alpha,q)=\cos{\lambda x}+\int^x_0 G(x,t)\cos{\lambda t}dt,
\end{equation}
where about the kernel $G(x,t)$ we know (in particular) that
\begin{equation}\label{eq13}
G(x,x)=-\cot\alpha + \frac{1}{2} \int^x_0 q(s)ds.
\end{equation}
It is also known that $G(x,t)$ satisfies to the Gelfand-Levitan integral equation
\begin{equation}\label{eq14}
G(x,t)+F(x,t)+\int^x_0 G(x,s)F(s,t)ds=0,\quad 0\leq t \leq x,
\end{equation}
where (see \cite{f1})
\begin{equation}\label{eq15}
F(x,t)=\sum_{n=0}^\infty
\Big( \frac{\cos{\lambda_n x}\cos{\lambda_n t}}{\tilde{a}_n}
-\frac{\cos{n x}\cos{n t}}{a_n^0}\Big)
\end{equation}
where $a_0^0=\pi$ and $a_n^0=\frac{\pi}{2}$ for $n=1,2,\dots$.
>From \eqref{eq13}--\eqref{eq15} it follows  that
\begin{equation}\label{eq16}
\begin{aligned}
G(0,0)&=-F(0,0)
=-\sum_{n=0}^\infty \Big( \frac{1}{\tilde{a}_n}-\frac{1}{a_n^0}\Big)\\
&=-\Big( \frac{1}{\tilde{a}_0}-\frac{1}{\pi}\Big)
- \sum_{n=1}^\infty \Big( \frac{1}{\tilde{a}_n}-\frac{2}{\pi}\Big)
 = -\cot\alpha.
\end{aligned}
\end{equation}
Thus, \eqref{eq10} is proved.

Let us now consider the functions (compare with \cite{j1})
\begin{equation}\label{eq17}
p(x,\mu_n)=\frac{\varphi(\pi-x,\mu_n,\alpha,q)}{\varphi(\pi, \mu_n,\alpha,q)}=
\frac{\varphi(\pi-x,\mu_n)}{\varphi(\pi, \mu_n)}, \quad n=0,1,2,\dots .
\end{equation}
Since $\varphi(x,\mu,\alpha,q)$ satisfies \eqref{eq1}, and
$$
p'(x,\mu_n)=-\frac{\varphi'(\pi-x,\mu_n)}{\varphi(\pi, \mu_n)}, \quad
p''(x,\mu_n)=\frac{\varphi''(\pi-x,\mu_n)}{\varphi(\pi, \mu_n)},
$$
we can see that $p(x,\mu_n)$ satisfies
$$
-p''(x,\mu_n)+q(\pi-x)p(x,\mu_n)=\mu_n p(x,\mu_n)
$$
and the initial conditions
\begin{equation}\label{eq18}
p(0,\mu_n)=1, \quad p'(0,\mu_n)=-\frac{\varphi'(\pi,\mu_n)}{\varphi(\pi,\mu_n)}
=-(-\cot\beta)=\cot\beta=-\cot(\pi-\beta).
\end{equation}
Also we have
\begin{gather*}
p(\pi,\mu_n)=\frac{\varphi(0,\mu_n)}{\varphi(\pi,\mu_n)}
=\frac{\sin\alpha}{\varphi(\pi,\mu_n)}
=\frac{\sin(\pi-\alpha)}{\varphi(\pi,\mu_n)},
\\
p'(\pi,\mu_n)=-\frac{\varphi'(0,\mu_n)}{\varphi(\pi,\mu_n)}
=-\frac{-\cos\alpha}{\varphi(\pi,\mu_n)}
=\frac{-\cos(\pi-\alpha)}{\varphi(\pi,\mu_n)}.
\end{gather*}
It follows, that $p_n(x):=p(x,\mu_n)$ satisfy to the boundary condition
$$
p_n(\pi)\cos(\pi-\alpha)+p'_n(\pi)\sin(\pi-\alpha)=0, \quad n=0,1,2,\dots.
$$
Let us denote $q^{*}(x):=q(\pi-x)$.
Since $\mu_n(q^{*},\pi-\beta,\pi-\alpha)=\mu_n(q,\alpha,\beta)$
(it is easy to prove and is well known \cite{i1}),
it follows, that $p_n(x)$, $n=0,1,2,\dots$, are the eigenfunctions of problem
$(q^{*},\pi-\beta,\pi-\alpha)$, which have the initial conditions \eqref{eq18}; i.e.
$p_n(x)=\tilde{\varphi}(x,\mu_n,\pi-\beta,q^{*})$, $n=0,1,2,\dots$.

Thus, as in \eqref{eq16}, for norming constants
$\hat{a}_n=\|p(\cdot, \mu_n)\|^2$ must satisfy
\begin{equation}\label{eq19}
\Big( \frac{1}{\hat{a}_0}-\frac{1}{\pi}\Big)
+ \sum_{n=1}^\infty \Big( \frac{1}{\hat{a}_n}-
\frac{2}{\pi}\Big) = \cot(\pi-\beta)=-\cot\beta.
\end{equation}
On the other hand, for the norming constants $\hat{a}_n$, according to
 \eqref{eq4}, \eqref{eq5} and \eqref{eq17}, we have
\begin{align*}
\hat{a}_n
&=\int_0^{\pi} p^2(x,\mu_n)dx\\
&= \int_0^{\pi} \frac{\varphi^2(\pi-x,\mu_n)}{\varphi^2(\pi,\mu_n)}dx\\
&=-\frac{1}{\varphi^2(\pi,\mu_n)} \int^0_{\pi} \varphi^2(s,\mu_n)ds\\
&= \frac{1}{\varphi^2(\pi,\mu_n)} \int^{\pi}_0 \varphi^2(s,\mu_n)ds\\
&=\frac{a_n(q,\alpha,\beta)}{\varphi^2(\pi,\mu_n)}
=\frac{\tilde{a}_n \sin^2\alpha}{c^2_n \sin^2\beta}=\tilde{b}_n.
\end{align*}
Therefore, we can rewrite \eqref{eq19} in the form
$$
\Big( \frac{1}{\tilde{b}_0}-\frac{1}{\pi} \Big)
- \sum_{n=1}^\infty \Big( \frac{1}{\tilde{b}_n}-\frac{2}{\pi}\Big)
 =\cot(\pi-\beta)= -\cot\beta.
$$
Thus, \eqref{eq11} holds, and Theorem \ref{thm2} is proved.
\smallskip

Let us note that the specification of the spectra
$\{\mu_n(q,\alpha,\beta) \} _{n=0} ^\infty $ (of a problem $(q,\alpha,\beta)$)
uniquely determines the characteristic function $\Phi(\mu)$
(see \cite[Lemma 2.2]{h1}, see also \cite[Lemma 1]{i1}), and 
its derivative $\frac{\partial \Phi(\mu)}{\partial \mu}=\dot{\Phi}(\mu)$
(see \cite[lemma2.3]{h1}).
In particular, if $\alpha,\beta \in (0,\pi)$ the following formulae hold:
\begin{gather}\label{eq20}
\dot{\Phi}(\mu_0)=-\pi \sin\alpha \sin\beta \prod_{k=1}^\infty
\frac{\mu_k - \mu_0}{k^2}, \\
\label{eq21}
\dot{\Phi}(\mu_n)=-\frac{\pi}{n^2} \left[ \mu_0 - \mu_n \right] \sin\alpha \sin\beta
\prod_{k=1, k \neq n}^\infty \frac{\mu_k - \mu_n}{k^2},
\end{gather}
for  $ n=1,2,\dots$.
On the other hand, it is easy to prove the relation
(see \cite[(2.16) in Lemma 2.2]{h1} and \cite[Lemma 1]{i1})
\begin{equation}\label{eq22}
a_n=-c_n  \dot{\Phi}(\mu_n).
\end{equation}

To take into account the relations \eqref{eq5} and \eqref{eq20}-\eqref{eq22}
 we find formulae for
$1/\tilde{b}_0$ and
$1/\tilde{b}_n$ with $n=1,2,\dots$
(in terms of $\{\mu_n\}_{n=0}^\infty$ and $\{\tilde{a}_n\}_{n=0}^\infty$):
\begin{gather}\label{eq23}
\frac{1}{\tilde{b}_0}=\frac{\tilde{a}_0}{\pi^2
 \big( \prod_{k=1}^\infty \frac{\mu_k - \mu_0}{k^2}\big)^2}, \\
\label{eq24}
\frac{1}{\tilde{b}_n}=\frac{\tilde{a}_n n^4}{\pi^2
 [\mu_0-\mu_n]^2  \big( \prod_{k=1, k \neq n}^\infty \frac{\mu_k - \mu_n}{k^2}\big)^2}.
\end{gather}
So, we can change the second assertion in Theorem \ref{thm2} by the assertion
\begin{align*}
&\frac{\tilde{a}_0}{\pi^2   \big( \prod_{k=1}^\infty \frac{\mu_k - \mu_n}{k^2}\big)^2}
-\frac{1}{\pi}\\
&+\sum_{n=1}^\infty\Big(
\frac{\tilde{a}_n n^4}{\pi^2  [\mu_0-\mu_n]^2
 \big( \prod_{k=1, k \neq n}^\infty \frac{\mu_k - \mu_n}{k^2}\big)^2}-
\frac{2}{\pi}\Big)=-\cot\beta,
\end{align*}
which coincides with \eqref{eq9}.

\section{proof of the Theorem \ref{thm1}}

For $\mu_n$ we have proved in \cite{h2}
(in a more general case, when $q\in L^1_{\mathbb{R}}[0, \pi]$)
the  asymptotic formula
\begin{equation}\label{eq25}
\mu_n(q,\alpha,\beta)=\left[ n+\delta_n(\alpha,\beta) \right] ^2 + \frac{1}{\pi}\int_0^{\pi} q(t)dt
+ r_n(q,\alpha,\beta),
\end{equation}
where $\delta_n$ is the solution of the equation
\begin{equation}\label{eq26}
\begin{aligned}
\delta_n(\alpha,\beta)
&=\frac{1}{\pi} \arccos{\frac{\cos\alpha}{\sqrt{\left[ n+\delta_n(\alpha,\beta)
 \right] ^2 \sin^2\alpha + \cos^2 \alpha}}}  \\
&\quad  -\frac{1}{\pi} \arccos{\frac{\cos\beta}{\sqrt{\left[ n+\delta_n(\alpha,\beta)
\right] ^2 \sin^2\beta + \cos^2 \beta}}}
\end{aligned}
\end{equation}
and $r_n(q,\alpha,\beta)=o(1)$, when $n \rightarrow \infty$, uniformly in
$\alpha, \beta \in [0,\pi]$ and $q$ from any bounded subset of
$L_{\mathbb{R}}^1[0,\pi]$ (we will write $q \in BL_{\mathbb{R}}^1[0,\pi]$).
It follows from \eqref{eq26} (see \cite{h2} for details), that if
$\sin \alpha \neq 0$ and $\sin \beta \neq 0$, $(\alpha,\beta \in (0,\pi))$, then
\begin{equation}\label{eq27}
\delta_n(\alpha,\beta)=\frac{\cot\beta - \cot\alpha}{\pi n} + O \big( \frac{1}{n^2} \big).
\end{equation}
It is not difficult to obtain from \eqref{eq25} that (see \cite{h3})
\begin{equation}\label{eq28}
\lambda_n=\sqrt{\mu_n}=n+\delta_n(\alpha,\beta)+\frac{\left[ q \right]}{2 \left[ n
+\delta_n(\alpha,\beta)\right]} + l_n + O \big( \frac{1}{n^2} \big),
\end{equation}
where
\[
l_n=\frac{1}{\pi [n+\delta_n(\alpha,\beta)]} \int_0^{\pi} q(x) \cos 2 \lambda_n x dx
=o \big( \frac{1}{n} \big)
\]
 and $[q]=\frac{1}{\pi}\int_0^{\pi} q(t) dt$.

In the case  $q\in L_{\mathbb{R}}^2[0, \pi]$ and $\alpha,\beta \in (0,\pi)$
it follows from \eqref{eq27} and \eqref{eq28} that
$l_n=\omega_n/n$, where $\left[ \omega_n \right] \in l^2$ and we can
rewrite \eqref{eq28} in the form
\begin{equation}\label{eq29}
\lambda_n=n+\frac{\omega}{n} +\frac{\omega_n}{n},
\end{equation}
where $\omega= \textrm{const} = \big(\cot\beta - \cot\alpha+ \frac{\pi}{2} [q]\big)/\pi$
and $\{\omega_n \} _{n=0}^\infty \in l^2$, i.e.
 $\sum_{n=1}^\infty | \omega_n |^2 < \infty$.
In \cite{f1} there is a proof of such assertion:

\begin{theorem}[\cite{f1}]\label{thm3}
For real numbers $\{\lambda_n^2 \} _{n=0}^\infty$ and $\{\tilde{a}_n \} _{n=0}^\infty$
to be the spectral data for a certain boundary-value problem $(q,\alpha,\beta)$
with $q\in L_{\mathbb{R}}^2[0, \pi]$, ($\alpha,\beta \in \left( 0,\pi \right)$), it is necessary
and sufficient that  relations \eqref{eq6} and \eqref{eq7}  hold.
\end{theorem}

Thus, if we have a real sequence
$\{\mu_n \}_{n=0}^\infty = \{\lambda_n^2 \}_{n=0}^\infty$,
which has the asymptotic representation \eqref{eq6}
and a positive sequence $\{\tilde{a}_n \}_{n=0}^\infty$,
which has the asymptotic representation \eqref{eq7},
then, according to the Theorem \ref{thm3}, there exist a function
$q\in L^2_{\mathbb{R}}[0, \pi]$ and some constants
$\tilde{\alpha},\tilde{\beta} \in (0,\pi)$ such that
$\lambda_n^2$, $n=0,1,2,\dots$, are the eigenvalues and
 $\tilde{a}_n$, $n=0,1,2,\dots$, are norming constants of a
Sturm-Liouville problem $(q,\tilde{\alpha},\tilde{\beta})$.

The function $q(x)$ and constants $\tilde{\alpha}, \tilde{\beta}$
are obtained on the way of solving the inverse problem by Gel'fand-Levitan method.
The algorithm of that method is as follows:

First we define the function $F(x,t)$ by formula \eqref{eq15}
 (note that this function is defined by $\{\lambda_n \}_{n=0}^\infty$ and
$\{\tilde{a}_n \}_{n=0}^\infty$ uniquely). Then we consider
the integral equation \eqref{eq14}, where $G(x,\cdot)$ is unknown function.
It is proved (see \cite{f1}) that provided \eqref{eq29} and \eqref{eq7}
the integral equation \eqref{eq14} has a unique solution $G(x,t)$.
With function $G(x,t)$, we construct a function
\begin{equation}\label{eq30}
\tilde{\varphi}(x,\lambda)=\cos{\lambda x}+\int^x_0 G(x,t)\cos{\lambda t}dt,
\end{equation}
which is defined for all $\lambda \in \mathbb{C}$. It is proved (see \cite{f1}) that
\begin{equation}\label{eq31}
- \tilde{\varphi}''(x,\lambda^2)+\Big( 2 \frac{d}{dx}G(x,x) \Big) \tilde{\varphi}(x,\lambda^2)
 = \lambda^2 \tilde{\varphi}(x,\lambda^2),
\end{equation}
almost everywhere on $(0,\pi)$,
\begin{gather*}
\tilde{\varphi}(0,\lambda^2)=1,\\
\tilde{\varphi}'(0,\lambda^2)=G(0,0).
\end{gather*}
If we state the condition
\begin{equation}\label{eq32}
G(0,0)=-\cot \alpha,
\end{equation}
then the solution \eqref{eq30} of  equation \eqref{eq31} will satisfy
the boundary condition \eqref{eq2}
$$
\tilde{\varphi}(0,\lambda^2) \cos \alpha+\tilde{\varphi}'(0,\lambda^2)\sin \alpha =0
$$
for all $\lambda \in \mathbb{C}$. Since from \eqref{eq14} it follows that
$G(0,0)=-F(0,0)$ and from \eqref{eq15}  that
$F(0,0)=-\sum_{n=0}^\infty \big( \frac{1}{\tilde{a}_n}-\frac{1}{a_n^0}\big)$,
we have that condition \eqref{eq32} can be represented as
$$
\sum_{n=0}^\infty \Big( \frac{1}{\tilde{a}_n}-\frac{1}{a_n^0}\Big)=\cot \alpha,
$$
which is our condition \eqref{eq8} on the sequence $\{\tilde{a}_n \}_{n=0}^\infty$.

It is also proved (see \cite{f1,z1}) that the expression
\[
 \frac {\tilde{\varphi}'_n (\pi)}{\tilde{\varphi}_n (\pi)}
= \frac {\tilde{\varphi}'(\pi, \lambda_n^2)}{\tilde{\varphi}(\pi,\lambda_n^2)}
\]
is a constant (i.e. does not depend on $n$),
 which we will denote by $-\cot \tilde{\beta}$.
So the functions $\tilde{\varphi}(x,\lambda_n^2), n=0,1,2,\dots$, are
the eigenfunctions of a problem $(q,\tilde{\alpha},\tilde{\beta})$, where
$q(x)=2 \frac{d}{dx} G(x,x)$, $\tilde{\alpha}$ is in advance given $\alpha$
and we want $\tilde{\beta}$ to be equals $\beta$.
We know from the Theorem \ref{thm2}, that for problem $(q,\alpha,\tilde{\beta})$
it holds
$$
\frac{1}{\tilde{b}_0}-\frac{1}{\pi}+\sum_{n=1}^\infty
\Big( \frac{1}{\tilde{b}_n}- \frac{2}{\pi} \Big)=- \cot \tilde{\beta}.
$$
Thus, if we obtain condition \eqref{eq11}, then we guarantee that
$\tilde{\beta}=\beta$. But \eqref{eq11} deals with the norming constants
$\tilde{b}_n$, which are not independent.
We have shown that we can represent $\tilde{b}_n$ by $\tilde{a}_n$ and
$\{\mu_k \}_{k=0}^\infty$
 (see the relations \eqref{eq23} and \eqref{eq24}).
 Therefore, instead of \eqref{eq11}, we obtain the condition in the
form \eqref{eq9}.
Theorem \ref{thm1} is proved.

\subsection*{Acknowledgments}
We would like to thank the anonymous referee
for pointing out that in an earlier version of Theorem \ref{thm1} we used
the condition $\tilde{a}_n=(\pi/2)+r_n$ with $r_n=O(1/n^2)$, from
an incorrect result  from
[Harutyunyan, T. N.; Asymptotics of the Norming Constants of 
the Sturm-Liouville Problem, Proceedings of YSU, 3, pp. 3-11, 2013].
The new version of Theorem \ref{thm1} uses
$\tilde{a}_n=(\pi/2)+(k_n/n)$ with $\{k_n\}\in \ell^2$,
which leads to the correct result.

This research is supported by the Open Society Foundations - Armenia,
 within the Education program, grant N18742.

\end{document}